\documentclass{svproc}
\usepackage{epsf}
\usepackage{graphicx}
\usepackage{amssymb}
\usepackage{hyperref}
\usepackage{todonotes}
\usepackage{amsmath}

\begin{document}
\def\l{\lambda}
\def\m{\mu}
\def\a{\alpha}
\def\b{\beta}
\def\g{\gamma}
\def\d{\delta}
\def\e{\epsilon}
\def\o{\omega}
\def\O{\Omega}
\def\v{\varphi}
\def\t{\theta}
\def\r{\rho}
\def\bs{$\blacksquare$}
\def\bp{\begin{proposition}}
\def\ep{\end{proposition}}
\def\bt{\begin{th}}
\def\et{\end{th}}
\def\be{\begin{equation}}
\def\ee{\end{equation}}
\def\bl{\begin{lemma}}
\def\el{\end{lemma}}
\def\bc{\begin{corollary}}
\def\ec{\end{corollary}}
\def\pr{\noindent{\bf Proof: }}
\def\note{\noindent{\bf Note. }}
\def\bd{\begin{definition}}
\def\ed{\end{definition}}
\def\C{{\mathbb C}}
\def\P{{\mathbb P}}
\def\Z{{\mathbb Z}}
\def\d{{\rm d}}
\def\deg{{\rm deg\,}}
\def\deg{{\rm deg\,}}
\def\arg{{\rm arg\,}}
\def\min{{\rm min\,}}
\def\max{{\rm max\,}}

% MATH -----------------------------------------------------------
\newcommand{\norm}[1]{\left\Vert#1\right\Vert}
\newcommand{\abs}[1]{\left\vert#1\right\vert}

\newcommand{\set}[1]{\left\{#1\right\}}
\newcommand{\setb}[2]{ \left\{#1 \ \Big| \ #2 \right\} }

\newcommand{\IP}[1]{\left<#1\right>}
\newcommand{\Bracket}[1]{\left[#1\right]}
\newcommand{\Soger}[1]{\left(#1\right)}

\newcommand{\Integer}{\mathbb{Z}}
\newcommand{\Rational}{\mathbb{Q}}
\newcommand{\Real}{\mathbb{R}}
\newcommand{\Complex}{\mathbb{C}}

\newcommand{\eps}{\varepsilon}
\newcommand{\To}{\longrightarrow}
\newcommand{\varchi}{\raisebox{2pt}{$\chi$}}

\newcommand{\E}{\mathbf{E}}
\newcommand{\Var}{\mathrm{var}}

% QED box --------------------------------------------------------
\def\squareforqed{\hbox{\rlap{$\sqcap$}$\sqcup$}}
\def\qed{\ifmmode\squareforqed\else{\unskip\nobreak\hfil
\penalty50\hskip1em\null\nobreak\hfil\squareforqed
\parfillskip=0pt\finalhyphendemerits=0\endgraf}\fi}

% This Document only ---------------------------------------------
% \renewcommand{\th}{^{\mathrm{th}}}
\newcommand{\Dif}{\mathrm{D_{if}}}
\newcommand{\Difp}{\mathrm{D^p_{if}}}
\newcommand{\GHF}{\mathrm{G_{HF}}}
\newcommand{\GHFP}{\mathrm{G^p_{HF}}}
\newcommand{\f}{\mathrm{f}}%fitting data
\newcommand{\fgh}{\mathrm{f_{gh}}}%fitting polynomial
\newcommand{\T}{\mathrm{T}}%Taylor polynomial vector
\newcommand{\K}{^\mathrm{K}}%Orde of the Taylor polynomial
\newcommand{\PghK}{\mathrm{P^K_{f_{gh}}}}%the generalized hermite interpolating polynomioal
\newcommand{\Dig}{\mathrm{D_{ig}}}%the digonal matrix with the singular values on the diagonal
\newcommand{\for}{\mathrm{for}}
\newcommand{\End}{\mathrm{end}}

\begin{titlepage}

\begin{center}

\topskip 5mm

{\LARGE{\bf {Lower bounds for high derivatives

\vskip 4mm

of smooth functions with given zeros}}}

\vskip 8mm

{\large {\bf G. Goldman}}

\vspace{6 mm}

{Department of Applied Mathematics, Tel Aviv University,
Tel Aviv 69978, Israel. e-mail: ggoldman@tauex.tau.ac.il}

\vspace{6 mm}

{\large {\bf Y. Yomdin}}

\vspace{6 mm}

{Department of Mathematics, The Weizmann Institute of Science,
Rehovot 76100, Israel. e-mail: yosef.yomdin@weizmann.ac.il}

\end{center}

\vspace{6 mm}
\begin{center}

{\bf Abstract}
\end{center}

{\small Let $f: B^n \rightarrow {\mathbb R}$ be a $d+1$ times continuously differentiable function on the unit ball $B^n$, with
$\max_{z\in B^n} |f(z)|=1$. A well-known fact is that if $f$ vanishes on a set $Z\subset B^n$ with a non-empty interior, then for each $k=1,\ldots,d+1$ the norm of the $k$-th derivative $\|f^{(k)}\|$ is at least $M=M(n,k)>0$. A natural question to ask is:

\medskip

{\it What happens for other sets $Z$? In particular, for finite, but sufficiently dense sets?}

\medskip

This question was partially answered in (\cite{Yom1},\cite{Yom5}-\cite{Yom7}). This study can be naturally related to a certain special settings of the classical Whitney's smooth extension problem.

\medskip

Our goal in the present paper is threefold: first, to provide an overview of the relevant questions and existing results in the general Whitney's problem. Second, we provide an overview of our specific setting and some available results. Third, we provide some new results in our direction. These new results extend the recent result of \cite{Yom6}, where an answer to the above question is given via the topological information on $Z$.}

\end{titlepage}

\newpage

%%%%%%%%%%%%%%%%%%%%%%%%%%%%%%%%%%%%%%%%%%%%%%%%%%%%%%%%%%%%%%%%%%%%%

\section{Introduction}\label{Sec:Intro}
\setcounter{equation}{0}

%One of the versions of the classical Whitney smooth extension problem is as follows: let $Z\subset B^n \subset {\mathbb R}^n$ be a closed subset of the unit ball $B^n$. For a given function $g:Z\to {\mathbb R}$ we look for $C^{d+1}$-smooth functions $f:B^n\to {\mathbb R}$, extending $g$ from $Z$ to $B^n$. In particular, if such functions $f$ exist (which is always the case for finite $Z$), we want to know the minimal possible $C^{d+1}$-norm of $f$.

\smallskip

%In dimension one, the classical work of Whitney \cite{Whi} provides a complete and explicit answer to this problem, in terms of divided finite differences, built on all the $(d+1)$-tuples of points in $Z$. Recent exciting developments in the general Whitney problem (see \cite{Fef} and references therein), provide a complete answer in any dimension. In particular, as in dimension one, it is enough to check only finite subsets of $Z$ with cardinality bounded in terms of $n$ and $d$ only. There is also an algorithmic way to estimate the minimal extension norm for any finite $Z$. However, a possibility of an explicit answer, as in dimension one, through a kind of multi-dimensional divided finite differences, remains an open problem.

\smallskip

In this paper we continue the study, started (chronologically) in (\cite{Yom1},\cite{Yom5}-\cite{Yom7}, \cite{Gol.Yom.10}, \cite{Gol.Yom.11}), of certain special settings of the classical Whitney's smooth extension problem (see \cite{Bru.Shv,Fef,Fef.Kla,Whi1,Whi2,Whi3}). 
%Our goal is twofold: first, to provide an overview of the relevant questions and existing results in general Whitney's problem, as well, as in our specific setting. Second, we provide some new results in our direction. \todo{The former is a repetition of the abstract}

\smallskip

The paper is organized as follows: in Section \ref{Sec:Whitney} we give (an extremely sketchy) overview of the Whitney's smooth extension problem.

\smallskip

In Section \ref{Sec:Rigidity} we provide (a slightly more extended) overview of our results in (\cite{Yom1},\cite{Yom5}-\cite{Yom7}, \cite{Gol.Yom.10}, \cite{Gol.Yom.11}), devoted to what we call ``smooth rigidity'', which is, essentially, a special case of the general Whitney's problem: we want to smoothly extend a zero function from a given set $Z$.

\smallskip

Finally, in Section \ref{Sec:New.Results} we present new results of this paper, which extend the results of \cite{Yom6}, bounding smooth rigidity in topological terms.

\section{Whitney's smooth extension problem}\label{Sec:Whitney}
\setcounter{equation}{0}

Let's recall shortly what this problem is about, and what is its current status.

\medskip

Let $E\subset B^n \subset {\mathbb R}^n$ be a closed subset of the unit ball $B^n$, and let $\bar f$ be a real function defined on $E$. Can $\bar f$ be extended to a $C^d$-smooth $f$ on ${\mathbb R}^n$, and, if the extension is possible, what is the minimal $C^d$-norm of $f$? (The $C^d$-norm of $f$ is accurately defined in Section \ref{Sec:Rigidity} below).

\medskip

Whitney himself, posing this problem, and partially answering it in \cite{Whi1,Whi2,Whi3}, provided two main results:

\medskip

\noindent 1. In dimension $n=1$, the necessary and sufficient extendability condition is (roughly) that for all the $d+1$ - tuples of points
$X=\{x_1,\ldots,x_{d+1}\} \in E$ the divided finite differences of $\bar f$ on $X$ are uniformly bounded (\cite{Whi2}).

\medskip

\noindent 2. In dimension $n\ge 2,$ Whitney considers ``more informative'' initial data: at each $x\in E$ the Taylor polynomial of degree $d$ of $f$ is given. This kind of data is usually called ``Whitney fields'' or ``Jet fields''). In \cite{Whi3} necessary and sufficient conditions are given for smooth extendability of the Whitney fields.

\medskip

We'll not touch Whitney fields anymore in this paper, instead we complete the discussion with the following two important remarks:

\medskip

\noindent A. The recent solution of the general Whitney problem in \cite{Fef} was, essentially, based on a reconstruction of the Whitney field of $\bar f$ from its point-wise value data.

\medskip

\noindent B. Whitney fields are important in many areas of mathematics. Let us just mention the profound work of Kolmogorov on the complexity of functional classes (see \cite{Kol.Tih}, and the references therein, and the huge research area opened by these works). For applications of Whitney fields in numerical analysis, see, for example, \cite{Wie.Yom}.

\medskip

Recent exciting developments in the general Whitney problem (see \cite{Bru.Shv,Fef,Fef.Kla} and references therein), provide essentially a complete answer to the general Whitney extension problem in any dimension. In particular, the results of \cite{Bru.Shv,Fef,Fef.Kla}), provide an important information on this problem, which was earlier available only in dimension one. The "finiteness principle" achieved in this recent work, claims that, as in classical Whitney's results in dimension one (\cite{Whi2}), it is enough to check only finite subsets of $Z$ with cardinality bounded only in terms of $n$ and $d$.

\medskip

For example, by \cite{Bru.Shv}, for $n=2, \ d=2$, i.e. for $2$-smooth functions of two variables, it is enough to check {\it the extendability of $\bar f$ from all the six-points subsets of $E$}. And the explicit criteria of the extendability from the six-points subsets are also given in \cite{Bru.Shv}.

\medskip

In general, an algorithmic way to provide the extension, and to estimate its minimal norm, is provided in \cite{Bru.Shv,Fef,Fef.Kla} and in other related publications.

\medskip

However, a possibility of an explicit answer, as in dimension one, through a kind of {\it multi-dimensional divided finite differences}, remains, in general, an open problem.

\section{Smooth rigidity}\label{Sec:Rigidity}
\setcounter{equation}{0}

Let us now turn to our work in (\cite{Yom1},\cite{Yom5}-\cite{Yom7}, \cite{Gol.Yom.10}). Our goal here was to consider some quite special cases of the general Whitney problem, but to provide more explicit answers. Specifically, let us describe in more detail the setting of the recent work \cite{Yom5}-\cite{Yom7}, which we use below. Let $Z\subset B^n \subset {\mathbb R}^n$ be a closed subset of the unit ball $B^n$. In \cite{Yom5}-\cite{Yom7} we look for $C^{d+1}$-smooth functions $f:B^n\to {\mathbb R}$, vanishing on $Z$. Such $C^{d+1}$-smooth (and even $C^\infty$) functions $f$ always exist, since any closed set $Z$ is the set of zeroes of a certain $C^\infty$-smooth function.

\medskip

We normalize the extensions $f$ requiring $\max_{B^n}|f|=1$. Put the pointwise norm $\|f^{(d+1)}(x)\|$ to be the sum of the absolute values of $d+1$-st order derivatives of $f$ at $x\in B^n$, and let the global norm $\|f^{(d+1)}\|$ be the maximum of $\|f^{(d+1)}(x)\|$ for all $x\in B^n$.

\medskip

We ask for the minimal possible norm of the last derivative $\|f^{(d+1)}\|$, which we call {\it the $d$-rigidity ${\cal RG}_d(Z)$ of $Z$}. In other words, for each normalized $C^{d+1}$-smooth function $f:B^n\to {\mathbb R},$ vanishing on $Z$, we have
$$
\|f^{(d+1)}\|\ge {\cal RG}_d(Z),
$$
and ${\cal RG}_d(Z)$ is the maximal number with this property.

\medskip

Bounding from below ${\cal RG}_d(Z)$ is, essentially, very close to the answering the Whitney problem for smooth extensions of the zero function $\bar f\equiv 0$ on $Z$.

\medskip

Papers \cite{Yom5,Yom6,Yom7}, provide certain bounds on ${\cal RG}_d(Z)$, in terms of the fractal geometry (or of the topology) of $Z$. Before stating these results more accurately, and in order to compare them with the general results, available today in Whitney's extension theory, let's make the following remark:

\medskip

Of course, the results of \cite{Fef,Fef.Kla} provide, in principle, an algorithmic way to estimate also our quantities ${\cal RG}_d(Z)$, for any closed $Z\subset B^n$ (via considering finite subsets of $Z$ of bounded cardinality). However, our goal in the present paper, as well as in our previous papers, related to Smooth rigidity, is somewhat different: {\it we look for an explicit answer, in terms of simple, and directly computable geometric (or topological) characteristics of $Z$}.

\medskip

Let's now come to specific results. As the ``model'' example, consider the case of dimension $n=1$. Here we have the following important fact, which can be proved by the standard interpolation (or finite difference) formulas:

\bp\label{prop:d.points}
For any $Z\subset B^1$ we have ${\cal RG}_d(Z)\ge \frac{(d+1)!}{2^{d+1}},$ if $Z$ consists of at least $d+1$ different points, and ${\cal RG}_d(Z)=0$ if $Z$ consists of at most $d$ different points.
\ep

%Thus in dimension one the minimal non-zero value of ${\cal RG}_d(Z)$ is $\frac{(d+1)!}{2^{d+1}}.$ This is not true any more in higher dimensions: for $Z\subset \hat B^n, \ n\ge 2,$ the $d$-rigidity ${\cal RG}_d(Z)$ attains arbitrarily small positive values (\cite{Yom6}. This fact is important for understanding the smooth rigidity phenomena, studied in the present paper, so in Section \ref{Sec:an.example} below we present in some detail an example.

%It is easy to see that there is {\it a uniform upper bound for the $d$-rigidity of all the subsets $Z\subset \hat B^n$.} Indeed, consider a certain $C^\infty$ function $\phi$, which vanishes identically on $\hat B^n$ and satisfies $M_0(\phi)=1$. Then $\psi$ vanishes in $Z$, and hence ${\cal RG}_d(Z)\le M_{d+1}(\phi).$

Another simple observation is the following:

\bp\label{prop:Z.interior}
For any $Z \subset B^n$ with a non-empty interior,
$$
{\cal RG}_d(Z) \ge \frac{(d+1)!}{2^{d+1}}.
$$
\ep
This fact easily follows from Proposition \ref{prop:d.points}. We just restrict any function $f$ to a certain straight line $\ell$, passing through $z_0$ with $|f(z_0)|=1,$ and through an interior point of $Z$.
%The present paper extend this result to all sufficiently dense $Z$.

\medskip

Let us mention here also an old result of \cite{Yom1}, related to smooth rigidity. Informally it can be stated as follows: if the set of zeros $Y(f)$ of a smooth function $f$ on $B^n$ does not look like a union of smooth hypersutfaces, of a total area bounded by a constant, depending only on $n,d$, then the norm of $f^{(d+1)}$ is not smaller than a certain positive constant, depending only on $n,d$.

\medskip

Many of the ``near-polynomiality'' results of \cite{Yom2,Yom.Com} can be naturally interpreted in terms of smooth rigidity. We plan to present some new results in this direction separately.

\subsection{Rigidity and Remez constant of $Z$}\label{Sec:Remez.zero.sets}

Here we present the results of \cite{Yom5}. We need a definition and some properties of the Remez (or Lebesgue, or norming, ...) constant (see, e.g. \cite{Rem,Bru.Yom,Yom3} and references therein).

\bd\label{Remez.constant}
For a set $Z\subset B^n \subset {\mathbb R}^n$ the Remez constant ${\cal R}_d(Z)$ is the minimal $K$
for which the inequality
$$
\sup_{B^n}\vert P \vert \leq K \sup_{Z}\vert P \vert
$$
is valid for any real polynomial $P(x)=P(x_1,\dots,x_n)$ of degree $d$.
\ed
%Thus for each polynomial $P$ of degree $d$ we have
%\be\label{eq:remez.ineq}
%\sup_{B^n}\vert P \vert \leq {\cal R}_d(Z) \sup_{Z}\vert P \vert.
%\ee
Clearly, we always have ${\cal R}_d(Z)\ge 1.$ For some $Z$ the Remez constant ${\cal R}_d(Z)$ may be equal to $\infty$. In fact, ${\cal R}_d(Z)$ is infinite if and only if $Z$ is contained in the set of zeroes
$$
Y_P=\{x\in {\mathbb R}^n, \ | \ P(x)=0\}
$$
of a certain polynomial $P$ of degree $d$. Sometimes it is convenient to use the inverse Remez constant $\hat {\cal R}_d(Z):=\frac{1}{{\cal R}_d(Z)}.$

%\subsection{Remez constant and rigidity (\cite{Yom5})}\label{Sec:main.results}
%\setcounter{equation}{0}

\medskip

In \cite{Yom4,Yom5} we show that the rigidity ${\cal RG}_d(Z)$ and the Remez constant $\hat {\cal R}_d(Z)$ are closely connected:

\begin{theorem}\label{thm:main1}(\cite{Yom5})
For any $Z \subset B^n$, \ \ \ $\frac{(d+1)!}{2}\hat {\cal R}_d(Z)\le {\cal RG}_d(Z)$.
\end{theorem}

The proof is heavily based on the result of \cite{Yom4}, published in the proceedings of the current series of conferences a few years ago.

\subsection{Rigidity and test curves}\label{Sec:test.curves}

As it was mentioned above, {\it for $Z$ with a non-empty interior we always have ${\cal RG}_d(Z)\ge \frac{(d+1)!}{2^{d+1}}$, independently of the size and the geometry of $Z$.}

\medskip

In \cite{Yom5,Yom6} the following question was discussed: {\it Can this last property be extended to other $Z$, beyond those with a non-empty interior? In particular, is it true for sufficiently dense finite sets $Z$?}

\medskip

A partial answer was given in \cite{Yom7}:

\begin{theorem}\label{thm:main.intro11}
If the box dimension $dim_e(Z)$ is greater than $n-\frac{1}{d+1}$, then
$$
{\cal RG}_d(Z)\ge M=M(n,d)>0,
$$
where the positive constant $M$ depends only on $n$ and $d$.
\end{theorem}
The box (or Minkowski, or entropy ...) dimension of a set $X$ is, informally, the power $\beta$ in the expression
$$
M(\e,X)\sim \left(\frac{1}{\e}\right)^\beta,
$$
where $M(\e,X)$ is the covering number of $X$ by the $\e$-balls.

\medskip

In particular, the result of Theorem \ref{thm:main.intro11} provides examples of discrete, but sufficiently dense, sets $Z$ for which ${\cal R}_d(Z)$ behaves in the same way as for sets with a non-empty interior.

\medskip

Let us describe our basic approach to the proof of this result in \cite{Yom7}. If we could find a straight line $\ell$ in ${\mathbb R}^n$, passing through the point $z_0$, where the absolute value $|f(z)|$ is equal to one, and through some $d+1$ distinct points in $Z$, we could immediately get the required lower bound ${\cal RG}_d(Z)\ge \frac{(d+1)!}{2^{d+1}}$ via the basic properties of $d$-rigidity, mentioned above (Proposition \ref{prop:d.points}).

\smallskip

However, for a generic finite set $Z$ any straight line $\ell$ meets $Z$ at one or two points at most. Instead we replace $\ell$ by a smooth curve $\o$, and try to mimic the calculations for $\ell$. This requires analysis of the high order chain-rule expressions, on one side, and construction of curves $\o$ with small high-order derivatives, passing through some $d+1$ distinct points in $Z$, on the other side.

\smallskip

We show, using a kind of ``discrete integral geometry'', that already for finite or discrete sets $Z$, which are dense enough (in particular, when $dim_e(Z) > n-\frac{1}{d+1}$), the required curves exist.

\medskip

We call the smooth curves $\o$ as above ``the test curves''. We expect such curves to play an important role in investigation of the Smooth Rigidity, the multidimensional divided finite differences (D.F.D's), and, ultimately, the general Whitney smooth extension problem.

\subsubsection{Polynomial test curves of higher degree}\label{Sec:testcurves.high.degree}

Another application of test curves is given in \cite{Gol.Yom.11} Here we assume that the set $Z$ of zeroes of $f$ sits on a polynomial (or near-polynomial) curve $\o$ of degree $s$, and provide certain specific smooth rigidity results.

\medskip

Our main results in \cite{Gol.Yom.11}, Theorem \ref{th:First.General.Ineq1} and Corollaries \ref{cor:basic1} and \ref{cor.example} below, can be considered as generalizations of Proposition \ref{prop:d.points} from straight lines to curves of higher degree.

\medskip

Consider polynomial parametric curves $\o$ of degree $s$ in ${\mathbb R}^n$. The curves $\o$ are given in the coordinate form by $\o(t)=(\o_1(t),\ldots,\o_n(t))$, with $t\in [-1,1]$, and with $\o_i(t)$ being polynomials in $t$ of the degree at most $s$. Denote, as usual, by $[\eta]$ the integer part of $\eta$.

\begin{theorem}\label{th:First.General.Ineq1}
Let $f: B^n \rightarrow {\mathbb R}$ and $\o$ be as above, with $\o([-1,1])\subset B^n$. Put $g(t)=f(\o(t))$. Then for each $t\in [-1,1]$ we have

\be\label{eq:main.ineq}
\sum^{d+1}_{|\alpha|=[\frac{d+1}{s}]+1} \|f^{(\alpha)}(\o(t))\|\ge C(n,d,s)\|g^{(d+1)}(t)\|.
\ee
with the positive constant $C(n,d,s)$, which is explicitly given in \cite{Gol.Yom.11}.

\medskip

In particular, we have

\be\label{eq:main.ineq1}
\sum^{d+1}_{|\alpha|=[\frac{d+1}{s}]+1} \|f^{(\alpha)}\|\ge C(n,d,s)\|g^{(d+1)}\|.
\ee
\end{theorem}

Thus for $\o$ a curve of degree $s>1$ we cannot translate the derivatives of $f(\o)$ directly to the derivatives of $f$ of the same order, as in the case of $\o$ being the straight line. But still we can bound from below the sum of the norms of the derivatives of $f$ of orders from $[\frac{d+1}{s}]+1$ to $d+1$ in terms of the derivatives of $g=f(\o)$.

\smallskip

As an immediate corollary we obtain:

\bc\label{cor:basic1}
In the assumptions as above, if the curve $\o$ crosses the zero set $Z$ of $f$ at at least $d+1$ points, and passes through a certain point $z_0$ with $|f(z_0)|\ge \gamma > 0$, then

\be\label{eq:main.ineq2}
\sum^{d+1}_{|\alpha|=[\frac{d+1}{s}]+1} ||f^{(\alpha)}||\ge C(n,d,s)\frac{\gamma(d+1)!}{2^{d+1}}.
\ee
\ec

%Let $\o$ be a polynomial parametric curve of degree $s$ in ${\mathbb R}^n$. So $\o$ is given by $\o(t)=(\o_1(t),\ldots,\o_n(t))$, with $\o_i(t)$ being polynomials in $t$ of the degree at most $s$. Assume that $\o([-1,1])\subset B^n$, but $\o([-1,1])$ is not contained in $Z$. Finally, we assume that $|Z\cap \o([-1,1])|=\infty$.

Another corollary is as follows:

\bc\label{cor.example}
Let $f: B^n \rightarrow {\mathbb R}$ be an infinitely differentiable function on $B^n,$ and let $\o$ be as above, with $\o([-1,1])\subset B^n$ not contained in the zero set $Z(f)$. Assume that $|Z\cap \o([-1,1])|=\infty$, i.e. there are infinitely many zeros of $f$ on the curve $\o$. Then there is an infinite number of the derivatives orders $m$ for which
$$
||f^{(m)}||\ge C(d,m,f,\o)\frac{(m+1)!}{2^{m+1}},
$$
with the positive constants $C(d,m,f,\o)$, which are explicitly given in \cite{Gol.Yom.11}.
\ec

\subsection{Rigidity via the topology of $Z$}\label{Sec:sing.top}

The following result is obtained in \cite{Yom6}: let $B^n$ be the unit $n$-dimensional ball. For a given integer $d$ let $Z\subset B^n$ be a smooth compact hypersurface with $N=(d-1)^n+1$ connected components $Z_j$ {\it with disjoint interiors $U_j$}. Let $\mu_j$ be the $n$-volume of the interior $U_j$ of $Z_j$, and put $\mu=\min \mu_j, \ j=1,\ldots, N$.

\begin{theorem}\label{thm:top.old}(\cite{Yom6})

Under the conditions above we have:

\medskip

\noindent 1. (Remez-type inequality). For each polynomial $P$ of degree $d$ on ${\mathbb R}^n$ we have
$$
\frac{\max_{B^n}|P|}{\max_{Z}|P|}\le (\frac{4n}{\mu})^d.
$$

\noindent 2. (Smooth rigidity). As a consequence, via \cite{Yom5}, we provide an explicit lower bound for the $(d+1)$-st derivatives of any smooth function $f$, which vanishes on $Z$, while being of order $1$ on $B^n$:
$$
||f^{(d+1)}||\ge \frac{1}{(d+1)!}(\frac{4n}{\mu})^d.
$$
\end{theorem}

\subsection{Rigidity via singular points and values}\label{Sec:sing.pts.val}

In a recent paper \cite{Gol.Yom.10} we bound from below higher derivatives of $f$ via the geometry of the critical points and values of $f$. While this setting does not exactly fit the rigidity definition, given above, it certainly remains in the general Whitney's problem framework. Basically, to conclude that a certain geometric configuration {\it of the critical points of $f$} implies lower bounds on the high-order derivatives of $f$, we apply in \cite{Gol.Yom.10} the results of \cite{Yom5,Yom6}, to the first order partial derivatives of $f$, instead of $f$ itself.

\medskip

In turn, to conclude that a certain geometric configuration {\it of the critical values of $f$} implies lower bounds on the higher derivatives of $f$, we try to ``read backward'' the (pretty old) results of \cite{Yom01,Yom.Com}. It turns out to be a non-trivial problem. We consider the results of \cite{Gol.Yom.10} as quite instructive, and plan to present further results in this direction separately.

\section{New results}\label{Sec:New.Results}
\setcounter{equation}{0}

Now we finally come to the new results of the present paper. These results form a direct continuation and extension of the ``topological'' results of \cite{Yom6}, presented in Section \ref{Sec:sing.top} above.

\medskip

Our topological assumption in \cite{Yom6} was that all the interiors of the components of the hypersurface $Z$ of zeroes of $f$ were disjoint. This assumption strongly simplifies the considerations, but also strongly restricts the applicability of the result. In the present paper {\it we drop this assumption, investigating arbitrary nesting configurations of the components of $Z$}. We do it here only on the plane, by two reasons: first, the presentation becomes truly elementary, and, second, pretty delicate topological considerations (in the line of the Jordan - Brower separation theorem in higher dimension) can be avoided. We strongly believe that our results below {\it can be verbally extended to higher dimensions}, and plan to present them separately.

\medskip

Let $Z=\cup^s_{i=1}Z_j$ be the union of the disjoint $s$ smooth ovals $Z_i$ in the plane ${\mathbb R}^2$. By ``oval'' we mean a closed curve without self-intersections in ${\mathbb R}^2$.

\medskip

For each two ovals $Z_{i_1}$ and $Z_{i_2}$ the important for us relations between them are the following:

\medskip

\noindent 1. $Z_{i_1}$ is entirely in the interior of $Z_{i_2}$ (or wise versa).

\medskip

\noindent 2. The interiors of $Z_{i_1}$ and of $Z_{i_2}$ are disjoint.

\medskip

Of course, here we strongly depend on the Jordan separation theorem, claiming that a continuous closed curve $S$ in ${\mathbb R}^2$ divides ${\mathbb R}^2$ into two parts - the exterior and the interior of $S$. However, with this theorem in dispose (which, of course, is the main topic, requiring a careful presentation in higher dimensions), we can provide a complete topological description of any collection of plane ovals.

\medskip

First we give a verbal topological description of $Z$: it consists of a certain number $s_0$ of disjoint ovals $O_l, \ l=1,\ldots, s_0$. In turn, each oval $O_l$ contains inside it ``depth one'' disjoint ovals $O_{l,l_.1}, \ l_1=1,\ldots, s_{l,1}$. We continue inductively, with the ``depth'' of the nested configuration of the ovals. See Fig. 1. We denote the total number of the ovals in $Z$ as $N$.

\medskip

Fig. 1. This is just an example of a nested configuration.

\medskip

Finally, we describe the compact domains $W_j$ whose boundaries are formed by the ovals $Z_i$. We start with the ``depth one'' ovals $O_l, \ l=1,\ldots, s_0$, and for each $O_l$ we consider the domain $W_l$, which is bounded by the oval $O_l$ from the outside, and by the depth two ovals $O_{l,l_1}, \ l_1=1,\ldots, s_{l,1}$ from the inside. Then we continue inductively. This construction defines for each $Z$ as above a unique collection ${\cal U}(Z)$ of the compact domains $W_j$, with non-intersecting interiors, whose boundaries are formed by the ovals $Z_i$. See Fig. 2.

\medskip

Fig. 2. This is a visual representation of the domains $W_j$, corresponding to the ovals configuration in Fig. 1. Notice the number of the domains is the same as the number of the ovals, confirming Lemma \ref{lem:top.iden} below in this special case.

\medskip

The fact important for our final result is the following:

\bl\label{lem:top.iden}
For each configuration $Z$ of plane ovals, the total number of the compact domains $W_j\in {\cal U}(Z)$, as above, is equal to the total number $N$ of the ovals in $Z$.
\el
\pr
We apply induction with respect to the nesting depth of $Z$. For all the ovals in $Z$ disjoint we are in situation of \cite{Yom6}, and the conclusion is immediate. Here we start the induction. As an example to the opposite situation, consider the ``totally nested case''. 

%shown in Fig.3.

%Finally, we describe the compact domains $W_j$ whose boundaries are formed by the ovals $Z_i$. We start with the ``depth one'' ovals $O_l, \ l=1,\ldots, s_0$, and for each $O_l$ we consider the domain $W_l$, which is bounded by the oval $O_l$ from the outside, and by the depth two ovals $O_{l,l_1}, \ l_1=1,\ldots, s_{l,1}$ from the inside. Then we continue inductively. This construction defines for each $Z$ as above a unique collection ${\cal U}(Z)$ of the compact domains $W_j$, with non-intersecting interiors, whose boundaries are formed by the ovals $Z_i$. See Fig. 2.

Now assume that the result is proved for a nesting depth $l$. To pass to a configuration of the nesting depth $l+1,$ we just add inside some ovals of the depth $l$ some new non-nested ovals (of the depth $l+1$). Next, in each oval of the depth $l$, we add its interior W', up to the new ovals, as one of the new domains $W_j$, and we add also all the interiors $W''$ of the new ovals of the nesting depth $l+1,$ (separately). 

\medskip

Clearly, if the number of the newly added ovals of the depth $l+1$ was $q$, the total number $N$ of the ovals increased by $q$. The total number of the new domains $W',W''$ is $q+1$, but one of these domains existed before the reconstruction. We conclude that the total number of the domains and of the ovals remains the same. This complete the proof. $\square$

\medskip

Let $\mu_j$ be the $n$-volume of the domain $W_j$. For $Z$ as above we put
$$
\mu(Z)=\min \mu_j, \  W_j \in {\cal U}(Z).
$$
Now we are ready to state and prove our main result:

\medskip

\begin{theorem}\label{thm:main.new}
For a given integer $d$ let $Z=\cup^s_{i=1}Z_j$ be the union of the disjoint $s\ge (d-1)^2+1$ smooth ovals $Z_i$ in the unit ball $B^2$ in the plane ${\mathbb R}^2$.

\medskip

Then (Remez-type inequality) for each polynomial $P$ of degree $d$ on ${\mathbb R}^2$ we have
$$
\frac{\max_{B^n}|P|}{\max_{Z}|P|}\le (\frac{8}{\mu(Z)})^d.
$$
As a consequence, via \cite{Yom5}, we have an explicit lower bound for the $(d+1)$-st derivatives of any smooth function $f$, which vanishes on $Z$, while being of order $1$ on $B^2$ (smooth rigidity):
$$
||f^{(d+1)}||\ge \frac{1}{(d+1)!}(\frac{8}{\mu(Z)})^d.
$$
\end{theorem}
\pr
We apply exactly the same arguments as in \cite{Yom6}, to the domains $W_j$. Thus, for the detailed proof we refer to \cite{Yom6}. The idea is that if the polynomial $P$, with $\max_{B^2} |P|$ is too small on $Z$, i.e. on the boundaries of the domains $W_j$, it still must be ``big'' inside each $W_j$. Otherwise it will be too small on the entire ball $B^2$, via the classical Brudny-Ganzburg-Remez inequality \cite{Bru.Gan}, applied separately to each $W_j$. Hence, being ``small'' on $Z$, but ``big'' inside each $W_j$, $P$ must have a critical point inside each $W_j$. As a result, we get more than $(d-1)^2$ critical point of $P$ - in contradiction with the Bezout theorem. $\square$

\medskip

Let us make a remark, omitted in \cite{Yom6}. The structure of the critical points of $P$ inside each $W_j$ may be rather complicated. In particular, they may form a variety $\Sigma$ of a positive dimension. But, by the construction, $\Sigma$ cannot cross the boundary of $W_j$, since the absolute value of $P$ on $\Sigma$ is a constant, strictly greater than the maximum of $|P|$ on the boundary of $W_j$the boundary of $W_j$.

\medskip

The number of such varieties $\Sigma$ can be explicitly bounded via the known results in semi-algebraic geometry. But we can stay with the Bezout bound of Theorem \ref{thm:main.new}, using the following standard trick: we make a small perturbation to $\tilde P=P+\xi$, adding to it an arbitrarily small linear form $\xi$, such that all the critical points of $\tilde P$ are non-degenerate. Since the perturbation $\xi$ is small, all the new non-degenerate critical points of $\tilde P$ remain inside the corresponding domains $W_j$. Now the usual Bezout bound can be applied.

\medskip

We conclude our discussion with the following open question: real algebraic geometry provides various restrictions on the topology of the sets of zeros of real polynomials of a given degree $d$. If for a given smooth $f$ any of these restriction is violated, then $f$ cannot be a polynomial of degree $d$, and hence the $d+1$-st derivative of $f$ is non-zero. The problem is {\it to provide a effective lower bound on this derivative}.

\medskip

As an example, the plane real curve $S$ of degree $6$ may have at most $11$ ovals, and the possible mutual positions of these ovals are given in the solution of the first part of Hilbert's 16-th problem. In particular, if the zero set $Z$ of $f$ contains $12$ or more ovals, then the $7$-st derivative of $f$ is non-zero. {\it But we are not aware of any approach to providing an effective lower bound.}

\medskip

In contrast, if there are $5^2+1=26$ ovals of the zero set $Z$ of $f$, then Theorem \ref{thm:main.new} provides an effective lower bound for the $7$-st derivative of $f$.

\medskip

Let us assume that zeros of $f$ are transversal, with an explicit lower bound $\gamma$ for the transversality measure. Then we expect that the results of \cite{Ler.Ste} would provide effective lower bounds for the higher derivatives of $f$, in case where the topology of $Z(f)$ violates the polynomial restrictions. However, these bounds blow up, as $\gamma \to 0$. Notice that in Theorem \ref{thm:main.new} there are no transversality assumptions.

%For a given integer $d$ let $Z\subset B^n$ be a smooth orientable compact hypersurface without boundary, with at least $(d-1)^n+1$ connected components. Under our assumptions each connected component $Z_i$ of $Z$, by itself, bounds exactly two domains in ${\mathbb R}^n$ - one compact, and one infinite.

%However, a given {\it collection of the connected components of $Z$} may bound several different compact domains - compare Fig. 1.

\medskip

%We consider all the possible compact domains $W_j \subset B^n$ with the non-empty interiors, whose boundaries are formed by some of the connected components of $Z$. Next we consider all the maximal nonintersecting unions (configurations) ${\cal U}$ of the domains $W_j$, containing
%$\nu({\cal U})$ domains, with $\nu({\cal U})\ge (d-1)^n+1$. ``Nonintersecting'' here means that for each $W_p, W_q \in {\cal U}$ we have $W_p \cap W_q = \emptyset$. ``Maximal'' means that no other $W_j$ can be added to ${\cal U}$, remaining nonintersecting. (There may be no such configurations at all). Compare Fig. 1.

\medskip

%We denote the set of the allowed configurations ${\cal U}$ by ${\cal V}(d,Z).$

%Let $\mu_j$ be the $n$-volume of the domain $W_j$. For each allowed configuration ${\cal U} \in {\cal V}(d,Z)$ we put
%$$\mu({\cal U})=\min \mu_j, \  W_j \in {\cal U}.$$ Finally, we define $\kappa_d(Z)$ as follows: \bd\label{def:kappa}
%For a given integer $d$ let $Z\subset B^n$ be a smooth compact hypersurface with at least $(d-1)^n+1$ connected components. Then the $d$-content $\kappa_d(Z)$ of $Z$ is defined as$$\kappa_d(Z) = \max_{{\cal U} \in {\cal V}(d,Z)} \ \mu({\cal U}).$$
%Here the maximum is taken over ${\cal U} \in {\cal V}(d,Z)$, i.e. over all the maximal nonintersecting configurations ${\cal U}$ of the domains $W_j$, containing
%\ed
\medskip

%\bt\label{thm:main.new}
%For a given integer $d$ let $Z\subset B^n$ be an oriented smooth compact hypersurface without boundary, and with at least $(d-1)^n+1$ connected components. Then (Remez-type inequality) for each polynomial $P$ of degree $d$ on ${\mathbb R}^n$ we have
%$$
%\frac{\max_{B^n}|P|}{\max_{Z}|P|}\le (\frac{4n}{\kappa_d(Z)})^d.
%$$
%As a consequence, via \cite{Yom5}, we have an explicit lower bound for the $(d+1)$-st derivatives of any smooth function $f$, which vanishes on $Z$, while being of order $1$ on $B^n$ (smooth rigidity):
%$$
%||f^{(d+1)}||\ge \frac{1}{(d+1)!}(\frac{4n}{\kappa_d(Z)})^d.
%$$
%\et
%\pr
%We apply exactly the same arguments as in \cite{Yom6}, to the domains $W_j$. Thus, for the detailed proof we refer to \cite{Yom6}. The idea is that if the polynomial $P$, with $\max_{B^n} |P|$ is too small on $Z$, i.e. on the boundaries of the domains $W_j$, it still must be ``big'' inside each $W_j$. Otherwise it will be too small on $B^n$, via the classical Brudny-Ganzburg-Remez inequality \cite{Bru.Gan}, applied separately to each $W_j$. Hence, being ``small'' on $Z$, but ``big'' inside each $W_j$, $P$ must have a critical point inside each $W_j$. As a result, we get more than $(d-1)^n$ critical point of $P$ - in contradiction with the Bezout theorem. $\square$

\medskip

%We provide now some comments, concerning the ``computability'' of the topological invariants $\nu({\cal U})$ and $\kappa_d(Z)$, introduced above.
%Let's start with the plane configurations. In this case the connected components of $Z$ are ovals $Z_l$ in ${\mathbb R}^2$. For each two ovals $Z_1$ and $Z_2$ the important for us relations between them are the following:

%On this base we can produce the purely combinatorial ``inclusion matrix'' $M(Z)$, presenting faithfully the topological relations between the ovals. Our claim, which we plan to carefully state and proof separately, is that the topological invariants $\nu({\cal U})$ and $\kappa_d(Z)$, are computable form the matrix $M(Z)$.

\medskip

%Moreover, for $n>2$ still the alternatives 1 and 2 above remain valid. Indeed, our connected components of $Z$ are oriented smooth compact hypersurfaces without boundary. Therefore, these hypersurfaces cannot cross one another, and the above alternative is preserved.

%and to the new results below. In the present paper we consider, as above, $C^{d+1}$-smooth functions $f:B^n\to {\mathbb R}$. But, in contrast to \cite{Yom5,Yom6,Yom7}, we do not consider zero sets of $f$. Instead, we assume some geometric conditions on the critical points and values of $f$, and, as above, derive in conclusion some lower bounds on $||f^{(d+1)}||$.

\medskip

\bibliographystyle{amsplain}

\end{document}